\documentclass[3p]{elsarticle}
\usepackage{amssymb, color}
\usepackage{amsmath}
 \usepackage{pifont}
 \usepackage[utf8]{inputenc}

\makeatletter
\def\LaTeX{\leavevmode L\raise.42ex
    \hbox{\kern-.3em\size{\sf@size}{0pt}\selectfont A}\kern-.15em\TeX}
\makeatother

\newcommand{\BibTeX}{{\rm B\kern-.05em{\sc
          i\kern-.025emb}\kern-.08em\TeX}}

\makeatletter
\def\@currentlabel{2.1}\label{e:dispaa}
\def\@currentlabel{2.21}\label{e:dispau}
\def\@currentlabel{2.22}\label{e:dispav}
\def\@currentlabel{2.23}\label{e:dispaw}
\def\@currentlabel{2.24}\label{e:dispax}
\def\theequation{\thesection.\@arabic\c@equation}
\makeatother

\renewcommand{\theequation}{\arabic{section}.\arabic{equation}}

\newcommand{\R}{\mathbb R}

\newcommand{\N}{\mathbb N}
\newcommand{\e}{\epsilon}

\def \D{\Delta}
\def \O{\Omega}
\newtheorem{thm}{Theorem} [section]
\newtheorem{lem}{Lemma} [section]
\newtheorem{prop}{Proposition} [section]

\newtheorem{definition}{Definition} [section]

\newtheorem{rmk}{Remark}[section]
\newtheorem{taggedtheoremx}{Theorem}
\newenvironment{taggedtheorem}[1]
 {\renewcommand\thetaggedtheoremx{#1}\taggedtheoremx}
 {\endtaggedtheoremx}

\renewcommand{\theequation}{\thesection.\arabic{equation}}
\renewcommand{\thesection}{\arabic{section}}
\renewcommand{\theequation}{\thesection.\arabic{equation}}
\let\ssection=\section\renewcommand{\section}{\setcounter{equation}{0}\ssection}

\def \p{\partial}

\begin{document}
\begin{frontmatter}

\title  { Solutions of super-linear elliptic equations and their Morse indices.}

\author[fd]{Foued Mtiri}
\ead{mtirifoued@yahoo.fr}

\address[fd]{IECL, UMR 7502, Universit\'{e} de Lorraine, France.}

\begin{abstract}

\noindent We investigate here the degenerate bi-harmonic equation:
$$\Delta_{m}^2 u=f(x,u)\; \;\;\mbox{in} \ \O,\quad u = \Delta u  = 0\quad \mbox{on }\; \p\Omega,$$ with  $m\ge 2,$  and also the degenerate
tri-harmonic equation: $$ -\Delta_{m}^3 u=f(x,u)\;\;\; \mbox{in} \ \O,\quad u = \frac{\p u}{\p \nu} =  \frac{\p^{2} u}{\p\nu^{2}} = 0\quad \mbox{on
}\; \p\Omega,$$ where $\Omega\subset \mathbb{R}^{N}$  is a bounded domain with smooth boundary $N>4$ or $N>6$ resp,  and $f \in
\mathrm{C}^{1}(\Omega\times \mathbb{R})$ satisfying suitable m-superlinear and subcritical growth conditions.  Our main purpose is to establish
$L^{p}$ and $L^{\infty}$ explicit bounds for weak solutions via the Morse index. Our results extend previous explicit estimate obtained in \cite{c,
HHF, hyf, lec}.
\end{abstract}
\noindent
\begin{keyword}
m-polyharmonic equation, Morse index, elliptic estimates.
\end{keyword}
\end{frontmatter}

\section{Introduction}

Consider the following m-polyharmonic equations

\begin{equation*}
(E_{k,m}): \quad \Delta_{m}^k u=f(x,u)\;\;\; \mbox{in}\,\, \ \O,
\end{equation*}
with the Dirichlet boundary conditions
\begin{align}
\label{D}
u = \frac{\p u}{\p \nu} = \ldots = \frac{\p^{k-1} u}{\p\nu^{k-1}} = 0  \quad \mbox{on }\; \p\Omega ,\ \ \ \ \mbox{if $k$ is odd};
\end{align}
or the Navier boundary conditions
\begin{eqnarray}\label{N}
u = \Delta u = \ldots = \D^{k-1} u = 0\quad \mbox{on }\; \p\Omega ,\ \ \ \ \mbox{if $k$ is even},
\end{eqnarray}
where $\Omega\subset \mathbb{R}^{N>2k}$ is a bounded domain with
smooth boundary. The m-poly-harmonic operator $\Delta_{m}^k $ is defined by
 \begin{align*}
\Delta_{m}^k u= \left\{\begin{array}{lll} D^{k} \left(|D^{k} u|^{m-2}D^{k} u\right),\ \ \ \ \  \ \ \mbox{if}\,\, k=2j,\\\\
\, -\mathrm{div}\left\{ \Delta^{\frac{k-1}{2}} \left(|D^{k} u|^{m-2}D^{k} u\right)\right\}\ \ \   \ \ \mbox{if}\,\, k=2j-1,
\end{array}
\right.
\end{align*}
where
\begin{equation*}
  D^k = \left\{\begin{array}{ll}\nabla\D^{j-1}& \mbox{if $k=2j-1$},\\
\D^{j}& \mbox{for $k=2j.$ }
\end{array}\right.
\end{equation*}
 The nonlinearity $f$ is a $\mathrm{C}^{1}(\Omega\times
\mathbb{R})$ function satisfying suitable superlinear and subcritical growth conditions. More precisely, we assume
\begin{enumerate}
  \item [$(H_1)$] (superlinearity) There exist $\mu> 0$  and $s_0> 0$ such that
  $$f'(x,s)s^{2}\geq (m-1+\mu)f(x,s)s > 0,\quad \mbox{for} \; |s|>s_{0}, \;\; x\in \Omega.$$
  \item [$(H_2)$] (subcritical growth) There exist $0<\theta<1$ and $s_0> 0$  such that
  $$\frac{mN}{N-km}F(x,s)\geq (1+\theta)f(x,s)s,\quad \mbox{for all} \; |s|>s_{0} \mbox{ and } x\in \Omega,$$
  where $F(x,s)=\displaystyle{\int_{0}^{t}}f(x,t)dt.$
  \item [$(H_3)$] There exist $C\geq 0$ and $s_0> 0$ such that
  $$|\nabla_x F(x,s)|\leq C(F(x,s)+1), \mbox{for all} \; |s|>s_{0} \mbox{ and }  x\in \Omega.$$
\end{enumerate}

We mention that under large growth conditions $(H_1)$-$(H_2)$ and if  we assume in addition that $f(x,s)$ grow less rapidly than $|s|^m$ near $0$ (respectively $f(x,.)$ is an odd function), then $(E_{k,m})$ has a nontrivial finite Morse index solution (respectively infinitely many finite Morse index solutions), obtained by mini-max method \cite{HH2}, (see also \cite{H2} for m=2).

\smallskip
For the second order case, i.e. $k = 1$ and $ m=2$, the equation $(E_{k,m}),$ becomes
\begin{equation*}
(E_{1,2}): \quad -\Delta u=f(x,u)\;\;\; \mbox{in}\,\, \ \O, \quad u= 0  \quad \mbox{on }\; \p\Omega.
\end{equation*}
 In \cite{c}, Bahri and Lions obtained the estimates of solutions in $H_0^1(\O)$ for superlinear and subcritical growth $f$, they used the blow-up
 technique and the boundedness of solutions' Morse index. Inspired by \cite{c}, Yang established in \cite{lec} the first explicit
 estimates of $L^{\infty}$ norm for solutions to $(E_{1,2})$  via the Morse index. Similar explicit estimates have been established in \cite{HHF} when the nonlinearity could be close to the critical growth. In particular the authors employed a cut-off function
with compact support to avoid the spherical integrals raised
in \cite{lec}  which are very difficult to control. The general higher order case is harder to achieve since we need to carefully handle some local interior estimate, especially near the boundary (see \cite{hyf} for the biharmonic and triharmonic cases under \eqref{N} and \eqref{D} with $k = 2$ and  $3$ respectively).

\smallskip
 However, when $\Omega$ is the entire space or the half space, Harrabi classified finite Morse index solutions of nonhomogeneous polyharmonic problem \cite{H} for all $k\geq 1$. His approach relies on a crucial idea, borrowed from \cite{RW}, where an appropriate family of test functions combined with an interpolation inequality related to weighted semi-norms are used to obtain the main integral estimate.  In contrast, explicit $L^\infty$-bounds in higher order do not seem to follow readily from similar arguments.

\medskip
It is a natural question  to ask if similar results can be observed for the  degenerate nonlinear operator $\Delta_{m}^k$ with $m\neq 2$. Very recently, in \cite{ka},  Hamdani and  Harrabi examined the case $k = 1$ and $m> 2$. They considered the following equation:

\begin{equation*}
(E_{1,m}): \quad -\Delta_{m} u=f(x,u)\;\;\; \mbox{in}\,\, \ \O, \quad u= 0  \quad \mbox{on }\; \p\Omega.
\end{equation*}
Under the above assumptions on the nonlinearity $f$ with $k=1$, they proved

 \begin{taggedtheorem}{A} Assume that f satisfies $(H_1)$-$(H_3)$ with $m> 2,$ then there exist positive
constant $C=(\Omega,f)$ such that any weak solution $u \in  C^{1,\alpha}_{loc}(\Omega)$ of $(E_{1,m})$ with finite Morse index $i(u)$, we have
$$\displaystyle{\int_{\Omega}}|\nabla u|^m dx\leq
C(i(u)+1)^{\alpha},\quad \|u\|_{L^{\infty}}\leq C (i(u)+1)^{\beta},$$
where
$$\alpha=\frac{m(\mu+m)}{\mu}+1 \quad \mbox{and}  \quad \beta= \dfrac{3m}{4
\theta(N-m)^{2}}\times \left(\frac{m(\mu+m)}{\mu}+1\right)\times\left(\frac{mN}{N-m}-1+q\right).$$
\end{taggedtheorem}

\medskip

In order to state our results more accurately, let us precie some basic definitions and notions.
Assume that $f$ satisfies the subcritical growth condition $(H_2).$
\begin{definition}\begin{itemize}\item The appropriate functional space of the variational setting of $(E_{k,m})$ is
\begin{equation*}\label{}
  \Sigma_k:= \left\{\begin{array}{llllllllll}
  W_0^{k,m}(\O):=\left\{v\in W^{k,m}(\O);\; \nabla^ju=0\;\mbox{on}\;\partial\Omega,\;\;\;\mbox{for}\;j=0, 1, . . ,k-1\right\}, & \mbox{if we work with \eqref{D}};\\\\
W_\vartheta^{k,m}(\O):=\left\{v \in W^{k,m}(\O), ;\;\Delta^jv=0\;\mbox{on}\;\partial\Omega,\;\mbox{for}\;\;\; j<\frac{k}{2}\right\}, & \mbox{if we work with \eqref{N}.}
\end{array}\right.
\end{equation*}
 \item We say that $u \in W^{k,m}(\Omega)$ is a weak solution of $(E_{k,m})$ if $u$ is a critical point of the following Euler-Lagrange energy functional
$$I(v)=\frac{1}{m}\int_{\O} |D^{k} v|^m dx-\int_{\O} F(x,v) dx, \quad \forall\; v \in \Sigma_k.$$

\item  For $m\geq 2$, we have $I\in C^2(\Sigma_k)$, and so the linearized operator of $(E_{k,m})$ at $u$ is given by
\begin{align*}
L_u(h,z)&:=\int_{\O} \Big[| D^{k} u|^{m-2} (D^{k} h\cdot D^{k} z) +(m-2)| D^{k} u|^{m-4} (D^{k} u\cdot D^{k} h)(D^{k} u\cdot D^{k}z)\Big]dx\\
& \; -\int_{\O}f'(x,u)hz dx  ,\; \; \; \forall (h,z) \in \Sigma^{2}_k,
\end{align*}
here $f'(x, u) :=\frac{\partial f}{\partial u}(x, u).$

\item  Let $u$ be a weak solution of $(E_{k,m}).$  The associated quadratic form to $(E_{k,m})$ of the linearized operator $L_u$ is defined by
\begin{align*}
  \Lambda_u (\phi)&:= \int_{\O}|D^{k} u|^{m-2}(D^k \phi)^2dx + (m-2)\int_{\O}|D^{k} u|^{m-4}(D^{k} u\cdot D^k \phi)^2dx \\
& \;-\int_{\O}f'(x,u)\phi^2dx\quad
  \mbox{for }\,\phi \in \Sigma_k,
\end{align*}

\item The Morse index of a classical solution $u$ of $(E_{k,m})$, denoted by $i(u)$
is defined as the maximal dimension of all subspaces of $\Sigma_k$ such
that $\Lambda_u(\phi) < 0$ in $\Sigma_{k}\setminus\{0\}.$ We say that $u$ is stable if its Morse index is equal to
zero.
\end{itemize}
\end{definition}

\begin{rmk}
\begin{itemize}
\item Observe that
\begin{align*}
  \int_{\O }|D^{k} u|^{m-4}(D^{k} u\cdot D^k \phi)^2 dx\leq \int_{\O }|D^{k} u|^{m-2} |D^k \phi|^2dx ,\; \forall \; \phi \in \Sigma_k.
\end{align*}
\item We should mention that when $\{u\phi_{j}\}_{1\leq j \leq i(u) + 1}$
are linearly independent, so there exists $j_{0}\in \{1,2,...,1+i(u)\}$ such that $\Lambda_u(u \phi_{j_0})\geq 0$
\begin{align}\label{quad}
 \int_{\O}f'(x,u)\phi^2dx\leq (m-1)\int_{\O}|D^{k} u|^{m-2}|D^k \phi|^2dx \quad \mbox{for }\,\phi \in \Sigma_k,
\end{align}
\item Observe that, besides the fact that many of our estimates work only in the case $m >2.$  Regarding in the case the case $1<m<2$, the energy functional $I$ belongs in $C^1(\Sigma_k)$ only, and in this case it is not clear which definition of stability would be the natural one.
\end{itemize}
\end{rmk}

\medskip
 Although we borrow many ideas from the previous works, we try to handle more general cases. In particular, we consider the degenerate bi-harmonic, i.e. when $k = 2,$ and also the degenerate  tri-harmonic (corresponding to
$k=3$) problems, under the Dirichlet or Navier boundary conditions, even if we believe that the results should hold true for more general $k \in\N.$

 \smallskip
From now on, we assume that $k = 2,$ or $3$ and  $m >2$. Our main objective is to obtain some $L^{p}$ and $L^{\infty}$ estimates for weak solutions of $(E_{2,m})$ and $(E_{3,m})$ via the Morse index. Our results read as follows

\begin{thm}\label{main2}
For $k = 2$ and  $3$ respectively,  there exists a positive constant $C=C(\Omega,\,f)$ such that if $u \in W^{2,m}(\Omega)$ is a weak solution of $(E_{2,m})$ with $m >2$ and $f \geq 0$ satisfying $(H_1)$-$(H_3)$ in $\R_+$; or if $u \in W^{3,m}(\Omega)$ is a weak solutions  of $(E_{3,m})$ with $f$ satisfying $(H_1)$-$(H_3)$, then

 $$\int_\Omega | D^k u|^m dx +\int_{\Omega}|f(x,u)|^{p_k} dx \leq C({i(u)}+1)^{\alpha_{k}},$$  where
$$p_{k,m}=\frac{mN}{N\left(m-(1+\theta)\right)+km(1+\theta))} \quad \mbox{and} \quad \alpha_{k,m}=\frac{km(2\mu+m)}{\mu}.$$
\end{thm}
By setting up a standard boot-strap iteration, as $f$ has subcritical growth, we can proceed similarly as in the proof of Theorem 2.2 in \cite{lec} and
claim that
\begin{thm}\label{main3}
If $u \in W^{2,m}(\Omega)\cap L^{\infty}(\Omega)$ is a weak solution of $(E_{2,m})$ with $m >2$ and  $f \geq 0$ satisfying $(H_1)$-$(H_3)$ in $\R_+$; or if
$u \in W^{3,m}(\Omega)\cap L^{\infty}(\Omega)$ is a weak solution of $(E_{3,m})$ with $f$ satisfying $(H_1)$-$(H_3)$, then there exists a
positive constant $C=C(\Omega,\,f)$ such that (for $k = 2$ or $3$ respectively),
$$ \|u\|_{L^{\infty}(\O)}\leq C (i(u)+1)^{\beta_{k,m}}, \quad \mbox{where } \; \beta_{k,m} =\frac{2k}{N}
\frac{\alpha_{k,m}}{p_{k,m}(2-p_{k,m})}\left[\frac{2k}{N(2-p_{k,m})}-\frac{1}{p_{k,m}}\right]^{-1},$$
where  $p_{k,m}$  and $\alpha_{k,m}$ is defined in Theorem \ref{main2}.
\end{thm}

Establishing some $L^{\infty}$ estimates which are \textbf{only} related to the Morse indices  for weak  solutions of $(E_{2,m})$ or $(E_{3,m})$  is
more
complicated, since these solutions are not $C^{2}(\Omega)$ or $C^{3}(\Omega)$ respectively. We shall derive a variant of the Pohozaev identity by using cut-off functions with compact
support. These functions allow us to avoid the spherical integral terms which appear in the classical Pohozaev identities and which are very difficult
to estimate, especially for   $(\D_{m})^2$ and  $(-\D_{m})^3$  situations.

\smallskip
 Moreover, We use  the quadratic form given by \eqref{quad} to
get some integral estimates, but the integration by parts argument yields  many terms which are difficult to control, even if we can borrow some ideas
from \cite{HHF, ka,hyf}, for example, the local $L^m$ norm of $\nabla u$ and $\nabla^2 u$, (see Lemma \ref{l.2.7ja} and \ref{l.2.7a} below). Even if
one proceeds similarly as for $(E_{2,m})$  or again  $(E_{3,m})$, there are additional difficulties that arise in each step.

 \smallskip
 Another difficulty, under $(H_1)$-$(H_3)$, the local $L^m$-estimate of $\nabla u$ and $\D u$ via the Morse index seem also hard to derive for
 $(E_{2,m})$ and $(E_{3,m})$. In fact we need to exhibit the explicit dependence on $i(u)$ (see Lemma \ref{l.2.3a} and \ref{l.2.3b} below).

\medskip

This paper is organized as follows: We give the proof of Theorem \ref{main2} for $k=2$ and $k=3$ respectively in sections 2 and 3. In the following,
$C$ denotes always a generic positive constant independent of the solution $u$, even  if their
value could be changed from one line to another one.

\section{Proof for  $k=2,$ and  $m> 2.$ }
\setcounter{equation}{0}

In order to prove our results, we need some technical lemmas which plays an important role in the proof of the above  theorem \ref{main2} for $k=2$. First,  remark that conditions  $(H_1)$ and $(H_2)$ imply that there exist two positive
constants $C_1$ and $C_2$ such that for $|s|$ large enough (resp. for $s$ large enough)
\begin{align}\label{B}
 \frac{(N-km)(1+\theta)}{mN} f(x,s)s-C_1 \leq F(x,s)\leq \frac{1}{m+\mu}f(x,s)s + C_1,
\end{align}
  \begin{align}\label{a}
 f(x,s)s \geq C_{1}(|s|^{m+\mu} - 1)
  \end{align}
  and
 \begin{align}\label{A}
  |f(x,s)|\leq C_{2}\left(|s|^{\frac{N\left(m-(1+\theta)\right)+km(1+\theta))}{(N-km)(1+\theta)}}+1\right).
\end{align}

Here we will prove Theorem \ref{main2} for $k = 2$.

\subsection{Preliminary technical results}

\smallskip

Let $y \in \mathbb{R}^N$ and $R>0.$ Throughout the paper, we denote by $B_R(y)$
the open ball of center $y $ and radius $R$ and $\partial\Omega_R(y):=\partial\Omega \cap B_R(y)$. For $x\in
B_{R}(y)\cap\Omega,$ let $n:=x-y.$ We denote also
$$u_{j_i\cdots j_k} :=\frac{\partial^{k} u}{\partial x_{j_1}\partial x_{j_2}\cdots\partial x_{j_k}}.$$

\medskip
First of all, we have the following Pohozaev identity.
\begin{lem}\label{l.2.1}
Let $u \in W^{2,m}(\Omega)$ be  a weak solution of $(E_{2,m})$ with $m >2$. Let $\psi \in C_c^{2}(B_R(y))$. Then
\begin{align*}
&\; \frac{Nm}{N-2m}\int_\O F(x,u)\psi dx +
\frac{m}{N-2m}\int_\O \nabla_{x} F(x,u)\cdot n\psi dx
-\int_\O (\Delta u)^{m}\psi dx \\
 =& \;-\frac{2m}{N-2m}\int_\O |\Delta u|^{m-2}\Delta u \nabla^{2}u (\nabla \psi, n) dx + \frac{1}{N-2m}\int_\O (\nabla\psi\cdot n)(\Delta u)^{m} dx\\
& \; -\frac{2m}{N-2m}\int_\O |\Delta u|^{m-2}\Delta u(\nabla u\cdot \nabla \psi) dx
- \frac{m}{N-2m}\int_\O |\Delta u|^{m-2}\Delta u(\nabla u\cdot n)\Delta \psi dx \\
& \; -\frac{m}{N-2m}\int_\O F(x,u)\nabla \psi\cdot n dx - \frac{m}{N-2m}\int_{\partial\Omega_R(y)} \frac{\p (|\Delta u|^{m-2}\Delta u)}{\p\nu}(\nabla u\cdot n)\psi d\sigma.
\end{align*}
\end{lem}

To describe our results more accurately, we need to make precise several terminologies. To establish a global estimate, we will cover the domain
 $\O$ by small balls and obtain local estimates.  To be more precise, consider $$\Omega_{1,R}:=
\left\{x \in \O: \mbox{ dist}(x,\partial \O)> \frac{R}{2}\right\}\;\;\mbox{and}\;\;\Omega_{2,R}:=\left\{x \in \O: \mbox{ dist}(x,\partial \O)\leq \frac{R}{3}\right\},\quad \forall\;  R>0. $$
The main difficulty is the estimates of $u$ near the boundary, that is, in $\Omega_{2,R}$. We need to choose carefully the balls as in \cite{lec}. Indeed, we will take balls with center lying in
\begin{align}\label{Gamma}
\Gamma(R):= \left\{x\in \mathbb{R}^N\backslash \O:\;\mbox{dist}(x,\partial \O)=
 \frac{R}{20}\right\}.
\end{align}
The domain $\O\backslash \Omega_{2,R}$ will be covered by balls with center lying in $\Omega_{1,R}.$ We can adapt the proof of Lemma 2.2 in \cite{hyf}  to obtain the following lemma  which is devoted to the control of the boundary term for $y \in \Gamma(R)$ in the above Pohozaev identity.
 \begin{lem}\label{l.2.2}
 Let $u \in W^{2,m}(\Omega)$ be  a weak solution of $(E_{2,m}),$ with $m >2.$  Assume that  $f(x,u)\geq 0,$ there exists $R_1(\O) > 0$ such that  for any $0< R \leq R_1(\O)$ and $y\in\Gamma(R),$ there holds
\begin{align}
 \nonumber\int_{\partial\Omega_R(y)}\frac{\p (|\Delta u|^{m-2}\Delta u)}{\p\nu}(\nabla u\cdot n)\psi d\sigma\geq0,
\end{align}
for any positive function $\psi \in C_c^{2}(B_R(y)).$
 \end{lem}

\smallskip

Consequently, we get
 \begin{lem}\label{c.2.1}
 Let $u \in W^{2,m}(\Omega)$ be  a weak solution of $(E_{2,m}),$ with $m >2$ and  $f \geq 0$ verifying $(H_{1})$-$(H_{3})$ in $\R_+$. Then for any $0< R\leq R_0$, $y \in\Gamma(R)$ and $0 \leq \psi \in C_c^{4}(B_R(y))$,  we conclude then
 \begin{align}
 \label{newl1}
 \begin{split}
& \int_\O f(x,u)u\psi dx + \int_\O (\D u)^m\psi dx \\
\leq& \; CR\|\nabla \psi\|_\infty \int_{A_{R,\psi}(y)}f(x,u) u dx + CR^{m}\int_{A_{R,\psi}(y)} |\nabla^2(u\nabla\psi)|^m dx\\
& \;  + C\Big(1 + R\|\nabla\psi\|_\infty\Big)\| \D u\|_{L^{m}(A_{R,\psi}(y))}^{m}+C\Big(R^m \|\nabla(\D\psi)\|_{\infty}^{m} + \|\Delta \psi\|_{\infty}^{m}\Big)\|u\|_{L^{m}(A_{R,\psi}(y))}^{m}\\
 & \; +CR^{m}\Big( \|\Delta \psi\|_{\infty}^m +\frac{ 1}{R^{m}}\|\nabla\psi\|_\infty^m+\|\nabla^{2}\psi\|_{\infty}^{m}\Big)\| \nabla u\|_{L^{m}(A_{R,\psi}(y))}^{m}  + + C(1 + \|\nabla \psi\|_\infty)R^{N+1},
\end{split}
\end{align}
where $$A_{R,\psi}(y) = B_{R}(y)\cap\Omega\cap\{\nabla \psi\neq0\}.$$
\end{lem}

\textbf{Proof.} Using Lemmas \ref{l.2.1}--\ref{l.2.2}, $(H_{1})$-$(H_{3})$ and \eqref{B}, we obtain

\begin{align}
\label{WXC}
\begin{split}
& (1+\theta)\int_\O f(x,u)u\psi dx - \int_\O (\D u)^m\psi dx \\
\leq & \;\frac{2m}{N-2m}
 \int_{A_{R,\psi}(y)}|\Delta u|^{m-1}|\nabla^{2}u (\nabla \psi, n)| dx + \frac{1}{N-2m}\int_{A_{R,\psi}(y)}(\Delta u)^{m}|\nabla\psi\cdot n| dx \\
& \; + \frac{2m}{N-2m}\int_{A_{R,\psi}(y)}|\Delta u|^{m-1}|\nabla u\cdot \nabla \psi| dx
+\frac{m}{N-2m}\int_{A_{R,\psi}(y)}|\Delta u|^{m-1}|\nabla u\cdot n||\Delta \psi| dx \\
& + \frac{1}{(N-2m)}\int_{A_{R,\psi}(y)}f(x,u)u|\nabla \psi\cdot n| dx + C R\int_{B_{R}(y)\cap\O}f(x,u)u \psi dx + + C(1 + \|\nabla \psi\|_\infty)R^{N+1}.
\end{split}
\end{align}
A direct calculation implies that
 \begin{align*}
\nabla^2u(\nabla\psi, n) = \sum_{ij} u_{ij}\psi_in_j = \sum_{ij}(u\psi_i)_{ij}n_j - u\nabla(\D\psi)\cdot n - \D\psi(\nabla u \cdot n) - \nabla^2\psi(\nabla u, n).
\end{align*}
By the Cauchy-Schwarz inequality, there exists $C > 0$ such that
\begin{align}
\label{KS}
\begin{split}
 \int_{A_{R,\psi}(y)}|\Delta u|^{m-1}| \nabla^{2}u (\nabla \psi, n)| dx
\leq& \; C\int_{A_{R,\psi}(y)}|\Delta u|^m dx + CR^{m}\int_{A_{R,\psi}(y)} u^{m}|\nabla(\D\psi)|^{m} dx \\
 &+ C R^{m}\int_{A_{R,\psi}(y)} |\nabla^2(u\nabla\psi)|^m dx\\
  & \; + CR^{m}\int_{A_{R,\psi}(y)}|\nabla u|^{m}\Big(\|\Delta \psi\|_\infty^m + \|\nabla^2\psi\|_\infty^m\Big) dx.
\end{split}
\end{align}
 Multiplying the equation  $(E_{2,m})$ by $u\psi,$ with $\psi \in C_c^4(B_R(y))$  and integrating by parts, we get readily
\begin{align}\label{SC}
\begin{split}
 \int_\O (\D u)^m\psi dx -\int_ \O f(x,u)u\psi dx &\leq C \int_{A_{R,\psi}(y)}|\Delta u|^{m-1}\Big[|\nabla u\cdot \nabla \psi| + |u||\Delta \psi| \Big] dx\\
 &\leq C\int_{A_{R,\psi}(y)}\Big[(\Delta u)^m + |\nabla u\cdot\nabla\psi|^m + (\Delta\psi)^mu^m\Big] dx.
 \end{split}
\end{align}
Remark that
\begin{align*}
\frac{\theta}{2} \int_\O (\D u)^m\psi dx + \frac{\theta}{2}\int_\O f(x,u)u\psi dx = & \; (1+\theta)\int_\O f(x,u)u\psi dx - \int_\O (\D u)^m\psi dx \\
& \; +\left(1 + \frac{\theta}{2}\right)\left[\int_{\O}(\D u)^m\psi dx - \int_{\O}f(x,u)u\psi dx\right].
\end{align*}
Fix $R_0 \in (0, R_1)$ such that $CR_0 < 1$.  Combining \eqref{WXC}-\eqref{SC}, using again Cauchy-Schwarz inequality, we get readily the estimate\eqref{newl1}. The proof for $y \in \O_{1, R}$ is completely similar, so we omit it.\qed

\medskip

To prove Theorem \ref{main2} for $k = 2$, we need also to establish an interior estimate. More precisely. Let $R>0,$ $y \in \Omega_{1,R}\cup \Gamma(R)$, $0<a<b$. Denote $$A:=A_{a}^{b}=\{x\in \mathbb{R}^{N};\; a<|x-y|<b\}, \quad A_\rho:= A_{a+\rho}^{b-\rho} \;\; \mbox{for }\; 0 < \rho < \frac{b-a}{4}.\eqno{(*)}$$

In the following Lemma, we establish an interior estimate for $\|\nabla u\|_{L^m(A_\rho\cap \O)},$ where we exhibit the dependence of the constant of this estimate with respect to $\rho.$

\begin{lem}\label{l.2.3a} There exists a constant
$C>0$ depending only on $N$ such that for any $u \in W^{2,m}(\O)\cap W_0^{1,m}(\O)$ and $0 < \rho < \min(1, \frac{b-a}{4}),$ we have
$$\|\nabla u\|_{L^m(A_\rho\cap \O)}^m \leq C\left(\frac{1}{\rho^{m}}\| u\|_{L^m(A\cap \O)}^m+\|\D u\|_{L^m(A\cap \O)}^m\right).$$
\end{lem}
\begin{rmk}\label{rem2.1}
If $f$ satisfies $(H_1)$ with $m >2$, using \eqref{a}, there holds
\begin{equation*}
\|u\|_{L^m(A\cap\O)}^m \leq C\left(\int_{A\cap \O}f(x,u)u dx\right)^{\frac{m}{m+\mu}}+ C.
 \end{equation*}
\end{rmk}
\subsection{Estimation via Morse index}
Let $u$  be a weak solution to $(E_{2,m})$ with finite Morse index $i(u)$. For $y \in \Gamma(R)\cup \O_{1, R}$, we denote
\begin{align}
 \label{defaj}
 A_j=:A_{a_j}^{b_j}\quad \mbox{with}  \;\;a_j = \frac{2(j+i(u))}{4(i(u)+1)}R,
\;\;b_j = \dfrac{2(j+i(u))+1}{4(i(u)+1)}R, \quad 1\leq j \leq i(u)+1.
\end{align}
Fix a cut-off function $\Phi \in C^\infty(\R)$ such that $\Phi =1$ in $[0, 1]$ and supp$(\Phi) \subset ({-\frac{1}{2}}, \frac{3}{2})$.
Let $$\phi_{j}(x):= \Phi\left(\dfrac{4 (i(u)+1)|x-y|}{R} - 2j - 2i(u)\right).$$
Then for any $1\leq j \leq i(u)+1$, $\phi_j \in C_c^\infty(B_R(y))$,
\begin{align}
\label{newest1}
\phi_{j}(x)=1  \;\; \mbox{in }\; A_{j},\quad \|\nabla\phi_{j}\|_{\infty}\leq \frac{C}{R}(1+i(u))\quad \mbox{and} \quad
\|\Delta\phi_{j}\|_{\infty}\leq \frac{C}{R^2} (1+i(u))^{2}.
\end{align}

\subsection{Main technical tool}

 As already mentioned, our proof of explicit estimates of $L^{p}$ and $L^{\infty}$  norm for weak solutions via the Morse index is based on the estimate $\Lambda_u (\phi)\geq0,$  with suitable test function. In fact, we have the following key  estimate which is an extension of the result in \cite{ka}.
  \begin{lem}
\label{lemnew} Let $u \in W^{2,m}(\O)$ be  a weak solution to $(E_{2,m}),$  with Morse index $i(u) < \infty.$ Assume that  $f$ satisfies $(H_1)$  with $m >2,$ then for any $0 < R \leq R_0$, $y\in \Gamma(R)\cup\O_{1, R}$, there exists a positive constant $C=C(\Omega,\,f)$ and  $j_{0}\in \{1,2,...,1+i(u)\}$ verifying
\begin{align}
\label{newest0}
\int_{A_{j_0}\cap\Omega}(\D
u)^{m} dx + \int_{A_{j_0}\cap\Omega}f(x,u)u dx \leq C
\left(\frac{1+i(u)}{R}\right)^\frac{2m(\mu+m)}{\mu}.
\end{align}
\end{lem}

 First, direct calculation shows that for $\epsilon \in (0, 1)$ and $\eta \in C_0^{\infty}(\O)$, with $0 \leq\eta \leq 1,$  there holds
\begin{align}
\label{newesTt47}
  \begin{split}
\int_{\Omega}|\Delta u|^{m-2} [\Delta (u\eta)]^2 dx & = \int_{\Omega}|\Delta u|^{m-2}\left(u \D \eta+2\nabla u\nabla \eta + \eta\D u\right)^{2} dx\\
& \leq \left(1+C\epsilon\right)\int_{\Omega}|\D u|^m {\eta}^2 dx + \frac{C}{\e}\int_{\Omega}|u|^2|\Delta u|^{m-2}
  |\D\eta|^2 dx \\
  &+\frac{C}{\epsilon} \int_{\Omega}|\Delta u|^{m-2} |\nabla u|^2|\nabla \eta|^2 dx.
    \end{split}
   \end{align}

   Take $\eta = \zeta^{2k}$ with $k > 2$, $0 \leq\zeta \leq 1$  and apply Young's inequality, we get
     \begin{align*}
  \int_{\Omega}|u|^{2}|\Delta u|^{m-2}|\D (\zeta^{2k})|^2 dx & \leq C_{k}\int_{\Omega}|u|^{2}|\Delta u|^{m-2}\left(|\D \zeta|^{2}+|\nabla \zeta|^4\right)\zeta^{\frac{4k(m-2)}{m}+\frac{2(4k-2m)}{m}} dx\\
& \leq C_{m, k}\epsilon^{2} \int_{\Omega}|\D u|^m \zeta^{4k} dx +
  C_{\epsilon, k, m}\int_{\Omega} |u|^m \left(|\D \zeta|^{2}+|\nabla \zeta|^4\right)^{\frac{m}{2}}\zeta^{4k-2m} dx,
\end{align*}
and
  \begin{align*}
  \int_{\Omega}|\Delta u|^{m-2} |\nabla u|^2|\nabla (\zeta^{2k})|^2 dx & =4 k^{2}\int_{\Omega}|\Delta u|^{m-2} |\nabla u|^2|\nabla \zeta|^2\zeta^{\frac{4k(m-2)}{m}+\frac{2(4k-m)}{m}} dx\\
& \leq  C_{m,k}\epsilon^{2} \int_{\Omega}|\D u|^m \zeta^{4k} dx +
  \frac{C_{m,k}}{\epsilon^{2}}\int_{\Omega} |\nabla u|^m |\nabla \zeta|^{m}\zeta^{4k-m} dx.
\end{align*}

Thus by the inequality \eqref{newesTt47}, together with these two estimates, one gets:
\begin{align}\label{newest4}
  \begin{split}
\int_{\Omega}|\Delta u|^{m-2} [\Delta (u\zeta^{2k})]^2 dx
  &\leq \left(1+C_{m, k}\epsilon\right) \int_{\Omega}|\D u|^m \zeta^{4k} dx+\frac{C_{m,k}}{\epsilon^{3}}\int_{\Omega} |\nabla u|^m |\nabla \zeta|^{m}\zeta^{4k-m} dx\\
  &+ C_{\epsilon, m, k}\int_{\Omega} |u|^m \left(|\D \zeta|^{2}+|\nabla \zeta|^4\right)^{\frac{m}{2}}\zeta^{4k-2m} dx.
  \end{split}
   \end{align}
   We will use also the following lemma
    \begin{lem}\label{l.2.7ja}
 Let $k \geq m/2 > 1$. For any $0<\epsilon < 1,$ there exists a positive constant $C >0$ such that for any $u\in W^{2,m}(\Omega)$ and $\zeta \in C_0^{\infty}(\O)$, with $0 \leq\zeta \leq 1,$   there holds
\begin{align*}
\begin{split}
 \int_{\Omega} |\nabla u|^m |\nabla \zeta|^{m}\zeta^{4k-m} dx
\leq  C\epsilon^{4}\int_{\Omega}| \D u|^{m}\zeta^{4k}dx
+ C \int_{\Omega}|u|^{m}\left(|\nabla \zeta|^{2m} +|\nabla^{2} \zeta|^{m}\right)\zeta^{4k-2m}dx.
\end{split}
\end{align*}
\end{lem}
\noindent{\bf Proof.} A simple calculation implies that
\begin{align*}
  \mathrm{div} \left(\nabla u|\nabla u|^{m-2}\right) u|\nabla \zeta|^{m}\zeta^{4k-m}= &\; \left(m-2\right)  \left(u |\nabla u|^{m-4}|\nabla \zeta|^{m}\nabla^{2}u (\nabla u, \nabla u)\zeta^{4k-m}\right)\\
 &\; +  u\D u|\nabla u|^{m-2}|\nabla \zeta|^{m}\zeta^{4k-m},
\end{align*}
and
\begin{align*}
u|\nabla u|^{m-2}\nabla u \cdot\nabla\left(|\nabla \zeta|^{m}\zeta^{4k-m}\right)=&\;  m u|\nabla u|^{m-2} |\nabla \zeta|^{m-2}\nabla^2\zeta(\nabla\zeta, \nabla u)\zeta^{4k-m} \\
+&(4k-m) u|\nabla u|^{m-2}|\nabla \zeta|^{m}(\nabla u\cdot \nabla\zeta)\zeta^{4k-m-1}.
\end{align*}
Hence, for any $0 \leq\zeta \leq 1,$   there exists $C > 0$ depending only on $k$ and $m$ such that
\begin{align}\label{KJHG}
\begin{split}
  A+B=&\;-\int_{\Omega} \mathrm{div} \left(\nabla u|\nabla u|^{m-2}\right) u|\nabla \zeta|^{m}\zeta^{4k-m}-\int_{\Omega} u|\nabla u|^{m-2}\nabla u \cdot\nabla\left(|\nabla \zeta|^{m}\zeta^{4k-m}\right) \\
 \leq &\;C_m\int_{\Omega}  |u|| \nabla^{2} u||\nabla u|^{m-2}|\nabla \zeta|^{m}\zeta^{4k-m}dx + \int_{\Omega}  |u||\D u||\nabla u|^{m-2}|\nabla \zeta|^{m}\zeta^{4k-m}dx\\
 &\;+ C_{ m,k} \int_{\Omega}|u||\nabla u|^{m-1}|\nabla \zeta|^{m-1}\left(|\nabla \zeta|^{2} +|\nabla^{2} \zeta|\right)\zeta^{4k-m-1}dx.
 \end{split}
\end{align}
For the last three  terms on the right hand side of \eqref{KJHG}. Applying Young's inequality, for any $\epsilon > 0$, there holds
 \begin{align*}
 \begin{split}
 & \;\int_{\Omega}  |u|| \nabla^{2} u||\nabla u|^{m-2}|\nabla \zeta|^{m-2+2}\zeta^{\frac{(4k-m)(m-2)}{m}+\frac{2(4k-m)}{m}}dx\\
 \leq &\;  C_{\epsilon, m}\int_{\Omega}  |u|^{\frac{m}{2}}| \nabla^{2} u|^{\frac{m}{2}}|\nabla \zeta|^{m}\zeta^{4k-m}dx + C\epsilon\int_{\Omega}|\nabla u|^{m}|\nabla \zeta|^{m}\zeta^{4k-m}dx\\
 \leq &\;C_{\epsilon, m}\int_{\Omega}  |u|^{m}|\nabla \zeta|^{2m}\zeta^{4k-2m}dx+ C\epsilon^{4}\int_{\Omega}| \nabla^{2} u|^{m}\zeta^{4k}dx+ C\epsilon\int_{\Omega}|\nabla u|^{m}|\nabla \zeta|^{m}\zeta^{4k-m}dx,
\end{split}
\end{align*}

\begin{align*}
 \begin{split}
 & \int_{\Omega}  |u|| \D u||\nabla u|^{m-2}|\nabla \zeta|^{m}\zeta^{\frac{(4k-m)(m-2)}{m}+\frac{2(4k-m)}{m}}dx\\
 \leq & \;C_{ \epsilon, m}\int_{\Omega}  |u|^{\frac{m}{2}}| \D u|^{\frac{m}{2}}|\nabla \zeta|^{m}\zeta^{4k-m}dx + C\epsilon\int_{\Omega}|\nabla u|^{m}|\nabla \zeta|^{m}\zeta^{4k-m}dx\\
 \leq &\;C_{\epsilon, m}\int_{\Omega}  |u|^{m}|\nabla \zeta|^{2m}\zeta^{4k-2m}dx+ C\epsilon^{4}\int_{\Omega}| \D u|^{m}\zeta^{4k}dx + C\epsilon\int_{\Omega}|\nabla u|^{m}|\nabla \zeta|^{m}\zeta^{4k-m}dx,
\end{split}
\end{align*}
and
 \begin{align*}
 \begin{split}
  &\;\int_{\Omega}|u||\nabla u|^{m-1}|\nabla \zeta|^{m-1}\left(|\nabla \zeta|^{2} +|\nabla^{2} \zeta|\right)\zeta^{\frac{(4k-m)(m-1)}{m}+\frac{(4k-2m)}{m}}dx\\
\leq &\; C_{\epsilon,m} \int_{\Omega}|u|^{m}\left(|\nabla \zeta|^{2m} +|\nabla^{2} \zeta|^{m}\right)\zeta^{4k-2m}dx
+ C\epsilon\int_{\Omega}|\nabla u|^{m}|\nabla \zeta|^{m}\zeta^{4k-m}dx.
 \end{split}
\end{align*}
Combining all these inequalities, we get the following inequality
\begin{align}\label{KJllHG}
\begin{split}
  A+B\leq &\;C_{\epsilon,m} \int_{\Omega}|u|^{m}\left(|\nabla \zeta|^{2m} +|\nabla^{2} \zeta|^{m}\right)\zeta^{4k-2m}dx
+ C\epsilon\int_{\Omega}|\nabla u|^{m}|\nabla \zeta|^{m}\zeta^{4k-m}dx\\
 &\;+ C\epsilon^{4}\int_{\Omega}| \D u|^{m}\zeta^{4k}dx+C\epsilon^{4}\int_{\Omega}| \nabla^{2} u|^{m}\zeta^{4k}dx.
 \end{split}
\end{align}

On the other hand, direct integrations by parts yield (recall that $u\in W^{2,m}(\Omega)$ )
\begin{align*}
  \begin{split}
  \int_{\Omega} |\nabla u|^m |\nabla \zeta|^{m}\zeta^{4k-m} dx = &\;\int_{\Omega} \nabla u\cdot\nabla u|\nabla u|^{m-2}|\nabla \zeta|^{m}\zeta^{4k-m}dx\\
 =& \; - \int_{\Omega}  \mathrm{div} \left(\nabla u|\nabla u|^{m-2}\right) u|\nabla \zeta|^{m}\zeta^{4k-m}dx\\
 &\; -\int_{\Omega} u|\nabla u|^{m-2}\nabla u \cdot\nabla\left(|\nabla \zeta|^{m}\zeta^{4k-m}\right)dx
 =: A + B.
 \end{split}
\end{align*}

By \eqref{KJllHG}, we deduce then
\begin{align}\label{KJllKKHG}
\begin{split}
&\; (1-C\epsilon) \int_{\Omega} |\nabla u|^m |\nabla \zeta|^{m}\zeta^{4k-m} dx\\
\leq &\;C\epsilon^{4}\int_{\Omega}| \nabla^{2} u|^{m}\zeta^{4k}dx+ C\epsilon^{4}\int_{\Omega}| \D u|^{m}\zeta^{4k}dx
+ C_{\epsilon,m} \int_{\Omega}|u|^{m}\left(|\nabla \zeta|^{2m} +|\nabla^{2} \zeta|^{m}\right)\zeta^{4k-2m}dx.
\end{split}
\end{align}

  Now we shall estimate the first term on the right hand side of \eqref{KJllKKHG}. Let $\psi^{r}\in C_0^{\infty}(\Omega),$ with  $r>2.$ By direct calculations, we get, as $u\in W^{2,m}(\Omega)$,

  \begin{align}\label{ecf}
|\nabla^2(u)|\psi^{r}
 \leq & \; C_{r}\Big[|u|\left(|\nabla \psi|^{2}\psi^{r-2}+|\nabla^{2} \psi|\psi^{r-1}\right)+ |\nabla u||\nabla \psi|\psi^{r-1}+ |\nabla^2(u\psi^{r})|\Big].
\end{align}
  Consider $\psi = \zeta,$ and  $r = \frac{4k}{m}\geq 2, $ so that $k\geq\frac{m}{2}.$ For any $0 \leq\zeta \leq 1,$   there exists $C_{m,k} > 0$ such that
 \begin{align}\label{0.2551l}
 \begin{split}
 &\int_{\Omega}| \nabla^{2} u|^{m}\zeta^{4k}dx\\
 \leq& \;C_{m,k}\int_{\Omega}\Big[|\nabla^2(u\zeta^{\frac{4k}{m}})|^m+|\nabla u|^{m}|\nabla \zeta|^{m}\zeta^{4k-m}+|u|^{m}\left(|\nabla \zeta|^{2m} +|\nabla^{2} \zeta|^{m}\right)\zeta^{4k-2m}\Big] dx.
 \end{split}
\end{align}
As $u\zeta \in W^{2, m}_0(\Omega)$, by standard elliptic theory, there exists $C_\O > 0$ depending only on $\Omega$ such that
 \begin{align}\label{0.j255l}
 \begin{split}
\int_{\Omega} |\nabla^2(u\zeta^{\frac{4k}{m}})|^m dx & \leq\;
C \int_{A_{j_0}\cap\Omega}|\D (u \zeta^{\frac{4k}{m}})|^m dx\\
\leq & \; C\int_{\Omega}\Big[|\D u|^{m}\zeta^{4k}+|\nabla u|^{m}|\nabla \zeta|^{m}\zeta^{4k-m}+|u|^{m}\left(|\nabla \zeta|^{2m} +|\nabla^{2} \zeta|^{m}\right)\zeta^{4k-2m}\Big] dx.
\end{split}
\end{align}

Combining \eqref{KJllKKHG}, \eqref{0.2551l} and \eqref{0.j255l}, we obtain
\begin{align*}
\begin{split}
 (1-C\epsilon) \int_{\Omega} |\nabla u|^m |\nabla \zeta|^{m}\zeta^{4k-m} dx
\leq  C\epsilon^{4}\int_{\Omega}| \D u|^{m}\zeta^{4k}dx
+ C \int_{\Omega}|u|^{m}\left(|\nabla \zeta|^{2m} +|\nabla^{2} \zeta|^{m}\right)\zeta^{4k-2m}dx.
\end{split}
\end{align*}
Take $\epsilon$ small enough, the lemma  follows.\qed

\smallskip
   Now, using Lemma \ref{l.2.7ja} , we obtain also
\begin{align}
\label{newest2}
\begin{split}
& \int_{\Omega}|\Delta u|^{m-2} [\Delta (u\zeta^{2k})]^2 dx\\
&\leq \; C\int_{\Omega} |u|^{m} \left(|\D \zeta|^{m}+|\nabla \zeta|^{2m} +|\nabla^{2} \zeta|^{m}\right)\zeta^{4k-2m} dx
+\left(1+C\epsilon\right) \int_{\Omega}|\D u|^m \zeta^{4k} dx .
 \end{split}
 \end{align}

Consider now the family of functions $\{u\phi_{j}^{k}\}_{1\leq j \leq i(u)+1}$, $k > 2$. With the definition of $\phi_j$, it's easy to see that different $\phi_j$ are supported by disjoint sets for different $j$, so they
are linearly independent as $u > 0$ in $\O$. Therefore, there must exist $j_{0}\in
\{1,2,...,1+i(u)\}$ such that $\Lambda_u(u \phi_{j_0}^{2k})\geq 0$ where $\Lambda$ is the quadratic form given by \eqref{quad}. Combining $\Lambda_u(u \phi_{j_0}^{2k} )\geq 0$ with \eqref{newest1} and  \eqref{newest2}, we obtain
    \begin{align}\label{0.2}
\displaystyle{\int_{\Omega}}f'
(x,u)u^{2}\phi_{j_{0}}^{4k} dx -\left(m-1\right)(1+C\epsilon)\int_{\Omega}|\D u|^m
\phi^{4k}_{j_{0}} dx \leq
\frac{C_\e}{R^{2m}}(1+i(u))^{2m}\int_{\Omega}|u|^{m}\phi_{j_{0}}^{4k-2m} dx.
\end{align}

Moreover, multiply the equation $(E_{2,m})$ by $u \zeta^{4k}$ and
integrate by parts, we get
\begin{align*}
 \int_{\Omega}|\D u|^{m} \zeta^{4k} dx-\int_{\Omega}  f(x,u)u\zeta^{4k} dx \leq\int_{\Omega}|u||\Delta u|^{m-1}
  |\D(\zeta^{4k})| dx
  + \int_{\Omega}|\Delta u|^{m-1} |\nabla u||\nabla (\zeta^{4k})| dx.
\end{align*}

Developing the right hand side, applying again  Lemma \ref{l.2.7ja}, we can conclude, for any $\epsilon>0$, there exists $C>0$
 such that
\begin{align}\label{f12}
\begin{split}
 \;\left(1- C\epsilon\right)  \int_{\Omega}|\D u|^{m} \zeta^{4k}-\int_{\Omega}  f(x,u)u\zeta^{4k} dx\leq
\frac{C_\e}{R^{2m}}(1+i(u))^{2m}\int_{\Omega}|u|^{m}\phi_{j_{0}}^{4k-2m} dx.
\end{split}
\end{align}

Now, take now $\zeta = \phi_{j_0}$  multiplying \eqref{f12} by $\frac{(m-1)(1+ 2C\epsilon)}{1- C\epsilon}$, adding it with \eqref{0.2} and $(H_1),$ we get
 \begin{align*}
& \;(m-1) C\epsilon\int_{\O}|\D u|^{m}\phi_{j_0}^{4k} dx + \left(\mu-\frac{C\epsilon(1-2m)-m+2}{1-C\epsilon}\right)\int_{\O}f(x,u)u\phi_{j_0}^{4k}dx\\
\leq& \;
\frac{C_\e}{R^{2m}}(1+i(u))^{2m}\int_{\Omega}|u|^{m}\phi_{j_{0}}^{4k-2m} dx +C_\e.
\end{align*}

Fix now $\epsilon <\frac{\mu-m+2}{C(\mu+2m-1})$, there holds

 \begin{align*}
\int_{\O}|\D u|^{m}\phi_{j_0}^{4k} dx +\int_{\O}f(x,u)u\phi_{j_0}^{4k}dx \leq \;
\frac{C_\e}{R^{2m}}(1+i(u))^{2m}\int_{\Omega}|u|^{m}\phi_{j_{0}}^{4k-2m} dx +C.
\end{align*}
Therefore, using \eqref{a} and $R \leq R_0$, for any $\epsilon'>0$,
  \begin{align}
  \label{ecf2}
 \begin{split}
\int_{\O}|\D u|^{m}\phi_{j_0}^{4k} dx+\int_{\O}f(x,u)u\phi_{j_0}^{4k}dx &\leq
C_{\epsilon'}\left(\frac{1+i(u)}{R}\right)^\frac{2m(\mu+m)}{\mu}+\e' \int_{\O}|u|^{\mu+m}\phi_{j_{0}}^{(4k-2m)(\frac{\mu+m}{m})}dx+ C \\
&\leq C_{\epsilon'}\left(\frac{1+i(u)}{R}\right)^\frac{2m(\mu+m)}{\mu}+ C\epsilon' \int_{\O}f(x,u)u\phi_{j_{0}}^{(4k-2m)(\frac{\mu+m}{m})}dx \\
& =C_{\epsilon'}\left(\frac{1+i(u)}{R}\right)^\frac{2m(\mu+m)}{\mu}+ C
\epsilon' \int_{\O}f(x,u)u\phi_{j_{0}}^{4k} dx.
\end{split}
\end{align}
For the last line, we used $(4k-2m)(\mu+m) = 4km$. Take $\epsilon' > 0$ small enough, the estimate \eqref{newest0} is proved. \qed

\subsection{Proof of Theorem \ref{main2} completed  }

\medskip
Now, we are in position to prove Theorem \ref{main2} for $k=2$. Fix
$$R = R_0, \quad \rho:=\frac{R}{10(i(u)+1)}, \quad A_{j_0,\rho}:=A_{a_{j_0}+\rho}^{b_{j_0}-\rho} \subset A_{j_0} \;\mbox{be as in $(*)$}.$$
According to Lemmas \ref{l.2.3a}, \ref{lemnew} and Remark \ref{rem2.1}, there exists a positive constant $C$ independent of $y \in \Gamma(R)\cup \O_{1, R}$ such that
 \begin{align}\label{3.7}
 \|\D u\|^m_{L^m(A_{j_0,\rho}\cap\O)}+ \|\nabla u\|^m_{L^m(A_{j_0,\rho}\cap\O)}\leq C
(1+i(u))^\frac{2m(\mu+m)}{\mu}.
\end{align}
Here, $a_{j_{0}}$ and $b_{j_{0}}$ are defined in \eqref{defaj} with $j_0$ given by Lemma \ref{lemnew}.

\medskip
Consider a cut-off function $\xi_{j_0} \in C_c^4(B_{b_{j_0}-\rho}(y))$ verifying
$\xi_{j_0}(x)\equiv 1$ in $B_{a_{j_0}+\rho}(y)$,
 with $$\|\nabla
\xi_{j_{0}}\|_{\infty}\leq \frac{C}{R}(1+i(u)), \quad \|\Delta \xi_{j_{0}}\|_{\infty}\leq
\frac{C}{R^{2}}(1+i(u))^2.$$
Applying Proposition \ref{c.2.1} with $\psi = \xi_{j_{0}}$, as $A_{R, \psi}(y) \subset A_{j_0, \rho}\cap \O$, we get
\begin{align}
\label{qssss}
\begin{split}
& \int_{\O}f(x,u)u\xi_{j_0} dx + \int_{\O}(\D u)^m\xi_{j_0} dx\\
 \leq & \;\ C(1+i(u))\int_{A_{j_0, \rho}\cap \Omega}\Big[(\D u)^m+f(x,u)u\Big] dx + C\int_{A_{j_0, \rho}\cap \Omega} |\nabla^2(u\nabla\xi_{j_{0}})|^m dx\\
  &+ C(1+i(u))^{3m} \| u\|_{L^{m}(A_{j_0, \rho}\cap\O)}^{m} + C(1+i(u))^{2m}\| \nabla u\|_{L^{m}(A_{j_0, \rho}\cap\O)}^{m}+ C(1+i(u))R^{N}.
\end{split}
\end{align}

Since  $u\nabla\xi_{j_0}=0$ on $\partial\O,$ by standard elliptic theory, there exists $C_\O > 0$ depending only on $\Omega$ such that
 \begin{align}\label{0.255}
 \begin{split}
\int_{\Omega} |\nabla^2(u\nabla\xi_{j_{0}})|^m dx & \leq C_\O \int_{A_{j_0, \rho}\cap\O}|\D (u \nabla\xi_{j_0})|^m dx\\
& \leq C\int_{A_{j_0, \rho}\cap\O} \Big[u^m |\nabla(\D\xi_{j_0})|^m + |\nabla u|^m|\nabla^{2} \xi_{j_0}|^m + (\D u)^m| \nabla\xi_{j_0}|^m\Big] dx.
\end{split}
\end{align}
From \eqref{qssss}, \eqref{0.255}, we get the following inequality
\begin{align}
\label{PL}
\begin{split}
& \int_{\O}f(x,u)u\xi_{j_0} dx + \int_{\O}(\D u)^m\xi_{j_0} dx\\
 \leq & \;\ C(1+i(u))\int_{A_{j_0, \rho}\cap \Omega}\Big[(\D u)^m+f(x,u)u\Big] dx+ C(1+i(u))^m\|\D u\|_{L^{m}(A_{j_0, \rho}\cap\O)}^{m}\\
  &+ C(1+i(u))^{3m} \| u\|_{L^{m}(A_{j_0, \rho}\cap\O)}^{m} + C(1+i(u))^{2m}\| \nabla u\|_{L^{m}(A_{j_0, \rho}\cap\O)}^{m}+ CR^{N}.
\end{split}
\end{align}
On the other hand, using Remark \ref{rem2.1} and Lemma \ref{lemnew}, there holds
\begin{align}\label{okm}
 \|u\|_{L^m(A_{j_{0}}\cap\O)}^m \leq C\left(\int_{A_{j_{0}}\cap \O}f(x,u)u dx\right)^{\frac{m}{m+\mu}} + C
 \leq C(1+i(u))^{\frac{2m^{2}}{\mu}}.
 \end{align}
Combining \eqref{newest0}, \eqref{3.7}, \eqref{PL} and  \eqref{okm}, one obtains
\begin{align*}
 \int_{\O}f(x,u)u\xi_{j_0} dx + \int_{\O}(\D u)^m\xi_{j_0} dx
\leq C(1+i(u))^\frac{2m(2\mu+m)}{\mu}.
\end{align*}
As $\frac{R}{2}< a_{j_{0}}$ and $R = R_0$, we get then for any $y \in \Gamma(R)\cup \O_{1, R}$,
\begin{align*}
\int_{B_{\frac{R_0}{2}}(y)\cap \O} \Big[|\D u|^{m} + f(x,u)u\Big] dx \leq C
(1+i(u))^\frac{2m(2\mu+m)}{\mu}.
\end{align*}
By covering argument and \eqref{A}, we get finally
$$\int_{\O} f(x, u)^{p_{2,m}} dx \leq C\int_\O f(x,u)u dx + C \leq C
(1+i(u))^{\alpha_{2,m}},$$
where $p_{2,m}=\frac{mN}{N\left(m-(1+\theta)\right)+2m(1+\theta))}$ and
$ \alpha_{2,m}=\frac{2m(2\mu+m)}{\mu}.$ So we are done.\qed

\smallskip\section{Proof of Theorem \ref{main2} for $m> 2,$ and $k=3$}
\setcounter{equation}{0} In this section, we consider the equation $(E_{3,m})$. We will proceed as for $(E_{2,m})$ and keep the same notations, but we replace the Navier boundary conditions by the Dirichlet boundary conditions and we have no more the sign condition for $f$.

\subsection{Preliminaries}

\medskip
 \begin{lem}\label{l.2.7a}
Let $k\geq\frac{3m}{4},$  and $\rho=R(1+i(u))^{-1}$. For any $0<\epsilon < 1$, there exists $C_{\e} >0$ such that for any $u\in W^{3,m}(\Omega)$ and $\zeta \in C_0^{\infty}(\O)$, with $0 \leq\zeta \leq 1,$  there holds
\begin{align*}
\begin{split}
&\rho^{-m}\int_{\O} |\nabla^{2} u|^m \zeta^{4k-m} dx +\rho^{-2m}\int_{\O} |\nabla u|^m \zeta^{4k-2m} dx+\rho^{-m}\int_{\O} |\D u|^m \zeta^{4k-m} dx\\
 \leq & \; C \epsilon^{3} \int_{\O}  |\nabla(\D u)|^m\zeta^{4k} dx + C_{\epsilon}\rho^{-3m} \int_{\O}|u|^{m}\zeta^{4k-3m}dx .
\end{split}
\end{align*}
\end{lem}
\noindent{\bf Proof .} We divide the proof in three parts.
\medskip

\noindent{\bf Step 1.}
 Using \eqref{ecf} with  $\psi = \zeta,$ and  $r = \frac{4k-m}{m}\geq 2, $ so that $k\geq\frac{3m}{4},$  we obtain
 \begin{align*}
 \begin{split}
 \int_{\Omega}| \nabla^{2} u|^{m}\zeta^{4k-m}dx\leq C\int_{\Omega}\Big[|\nabla^2(u\zeta^{\frac{4k-m}{m}})|^m+|\nabla u|^{m}|\nabla \zeta|^{m}\zeta^{4k-2m}+|u|^{m}\left(|\nabla \zeta|^{2m} +|\nabla^{2} \zeta|^{m}\right)\zeta^{4k-3m} \Big]dx.
 \end{split}
\end{align*}
 Since $u\zeta \in W^{3, m}_0(\Omega)$, there exists $C>0$ depending only on $\O$ such that
 \begin{align*}
 \begin{split}
 &\;\int_{\Omega}| \nabla^{2} u|^{m}\zeta^{4k-m}dx\\
 &\leq\; C\int_{\Omega}\Big[|\nabla^2(u\zeta^{\frac{4k-m}{m}})|^m+|\nabla u|^{m}|\nabla \zeta|^{m}\zeta^{4k-2m}+|u|^{m}\left(|\nabla \zeta|^{2m} +|\nabla^{2} \zeta|^{m}\right)\zeta^{4k-3m} \Big]dx\\
 &\leq\;C \int_{A_{j_0}\cap\Omega}|\D(u\zeta^{\frac{4k-m}{m}})|^m+ C\int_{\Omega}\Big[|\nabla u|^{m}|\nabla \zeta|^{m}\zeta^{4k-2m}+|u|^{m}\left(|\nabla \zeta|^{2m} +|\nabla^{2} \zeta|^{m}\right)\zeta^{4k-3m} \Big]dx\\
 &\leq\; C\int_{\Omega}\Big[|\D u|^{m}\zeta^{4k-m}+|\nabla u|^{m}|\nabla \zeta|^{m}\zeta^{4k-2m}+|u|^{m}\left(|\nabla \zeta|^{2m} +|\nabla^{2} \zeta|^{m}\right)\zeta^{4k-3m} \Big]dx.
 \end{split}
\end{align*}
So we get
 \begin{align}\label{0.2551}
\rho^{-m}\int_{\O} |\nabla^{2} u|^m \zeta^{4k-m} dx\leq\; C\int_{\Omega}\Big[\rho^{-m}|\D u|^{m}\zeta^{4k-m}+\rho^{-2m}|\nabla u|^{m}|\zeta^{4k-2m}+\rho^{-3m}|u|^{m}\zeta^{4k-3m} \Big]dx.
\end{align}
\noindent{\bf Step 2.}
A simple calculation implies that
\begin{align}\label{0.2551234}
\begin{split}
 &\rho^{-2m}\int_{\O} |\nabla u|^m \zeta^{4k-2m} dx\\
 & \;= - \rho^{-2m}\Big[\int_{\O}  u\;  \mathrm{div} \left(\nabla u|\nabla u|^{m-2}\right)\zeta^{4k-2m}dx+ \int_{\O} u|\nabla u|^{m-2}\nabla u \cdot\nabla\left(\zeta^{4k-2m}\right)dx\Big].
 \end{split}
\end{align}
hence the first term on the right hand side of \eqref{0.2551234} can be estimated as
 \begin{align*}
 &- \rho^{-2m}\int_{\O}  u\;  \mathrm{div} \left(\nabla u|\nabla u|^{m-2}\right)\zeta^{4k-2m} dx\\
 \leq& \; C_{m}\rho^{-2m}\int_{\O}  |u|| \nabla^{2} u||\nabla u|^{m-2}\zeta^{4k-2m}dx+\rho^{-2m}\int_{\O}  |u||\D u||\nabla u|^{m-2}\zeta^{4k-2m}dx.
\end{align*}
Applying Young's inequality, for any $\epsilon > 0$, there holds
 \begin{align}\label{0.25512vc}
 \begin{split}
 &\rho^{-2m}\int_{\O}  |u|| \D u||\nabla u|^{m-2}\zeta^{4k-2m}dx= \;\rho^{-2(m-2)-4}\int_{\O}  |u|| \D u||\nabla u|^{m-2}\zeta^{\frac{(4k-2m)(m-2)}{m}+\frac{2(4k-2m)}{m}}dx\\
 \leq &\;C _{\epsilon}\rho^{-3m}\int_{\O}  |u|^{m}\zeta^{4k-3m}dx+ C \epsilon\rho^{-m}\int_{\O}| \D u|^{m}\zeta^{4k-m}dx+ C \epsilon\rho^{-2m}\int_{\O}|\nabla u|^{m}\zeta^{4k-2m}dx.
\end{split}
\end{align}
and
 \begin{align}\label{0BN.25512}
 \begin{split}
 & \;C_{m}\rho^{-2m}\int_{\O}  |u|| \nabla^{2} u||\nabla u|^{m-2}\zeta^{\frac{(4k-2m)(m-2)}{m}+\frac{2(4k-2m)}{m}}dx\\
 \leq &\;C_{\epsilon}\rho^{-3m}\int_{\O}  |u|^{m}\zeta^{4k-3m}dx+C\epsilon\rho^{-m}\int_{\O}| \nabla^{2} u|^{m}\zeta^{4k-m}dx+ C\epsilon\rho^{-2m}\int_{\O}|\nabla u|^{m}\zeta^{4k-2m}dx.
\end{split}
\end{align}
Combining \eqref{0.2551}, \eqref{0.25512vc} and \eqref{0BN.25512}, we obtain the estimate for the first left term in \eqref{0.2551234}:
 \begin{align}\label{SSS0.255}
 \begin{split}
&\;-\rho^{-2m} \int_{\O}  u  \;\mathrm{div} \left(\nabla u|\nabla u|^{m-2}\right)\zeta^{4k-m} dx \\
 \leq &\; C\epsilon\rho^{-m}\int_{\O}|\D u|^{m}\zeta^{4k-m} dx + C\epsilon\rho^{-2m}\int_{\O}|\nabla u|^{m}\zeta^{4k-2m}dx+C_{\epsilon}\rho^{-3m} \int_{\O}|u|^{m}\zeta^{4k-3m}dx.
 \end{split}
\end{align}

On the other hand, by  Young’s inequality, we get, for any $\epsilon >0$
 \begin{align}\label{SS0.255}
 \begin{split}
&-\rho^{-2m} \int_{\O} u|\nabla u|^{m-2}\nabla u \cdot\nabla\left(\zeta^{4k-2m}\right)dx\\
 \leq &\;C\rho^{-3m} \int_{\O}|u|^{m}\zeta^{4k-3m}dx+C\epsilon\rho^{-2m}\int_{\O}|\nabla u|^{m}\zeta^{4k-2m}dx.
 \end{split}
\end{align}
Combining \eqref{SSS0.255}--\eqref{SS0.255}, one obtains
  \begin{align}\label{WE.0.0a}
\rho^{-2m}(1 - C\epsilon)\int_{\O} |\nabla u|^m \zeta^{4k-2m} dx
 \leq C_{\epsilon}\rho^{-3m} \int_{\O}|u|^{m}\zeta^{4k-3m}dx +C \epsilon \rho^{-m}\int_{\mathbb{R}^N}| \D u|^{m}\zeta^{4k-m}dx.
\end{align}

\noindent{\bf Step 3.}
 On the other hand, direct integrations by parts yield
\begin{align}
\label{11.0.0a}
\begin{split}
&\rho^{-m} \int_{\O}|\Delta u|^{m}\zeta^{4k-m} dx \\
= &\; -\rho^{-m}\Big[(m-1) \int_{\O}|\Delta u|^{m-2}\nabla u\cdot\nabla(\D u)\zeta^{4k-m} dx
 + (4k-m)\int_{\O}|\Delta u|^{m-1} \nabla u\cdot\nabla \zeta \zeta^{4k-m-1} dx\Big]\\
\leq &\; C\rho^{-m} \int_{\O}|\Delta u|^{m-2}|\nabla u||\nabla(\D u)|\zeta^{\frac{(4k-m)(m-2)}{m}+\frac{2(4k-m)}{m}} dx \\
+&\;C\rho^{-(m-1)-2}\int_{\O}|\Delta u|^{m-1} |\nabla u|\zeta^{\frac{(4k-m)(m-1)}{m}+\frac{(4k-2m)}{m}} dx.
\end{split}
\end{align}
Applying Young’s inequality, we get, for any $\epsilon >0$
\begin{align*}
\rho^{-m}(1-C \epsilon)\int_{\O}|\Delta u|^{m}\zeta^{4k-m} dx
\leq & \; \frac{C}{\epsilon^{3}}\Big[\rho^{-m} \int_{\O}|\nabla (\D u)|^{\frac{m}{2}} |\nabla u|^{\frac{m}{2}}\zeta^{4k-m} dx +\rho^{-2m} \int_{\O} |\nabla u|^{m}\zeta^{4k-2m} dx\Big]\\
\leq & \;C \epsilon^{3} \int_{\O}  |\nabla(\D u)|^m\zeta^{4k} dx  + C_{\epsilon}\rho^{-2m} \int_{\O} |\nabla u|^{m}\zeta^{4k-2m}dx.
\end{align*}

Take another small enough $\e$ in \eqref{WE.0.0a}, there holds

\begin{align*} \rho^{-m}(1-C \epsilon)\int_{\O} |\D u|^m \zeta^{4k-m} dx \leq C \epsilon^{3} \int_{\O}  |\nabla(\D u)|^m\zeta^{4k} dx + C_{\epsilon}\rho^{-3m} \int_{\O}|u|^{m}\zeta^{4k-3m}dx .
\end{align*}
The proof is completed. \qed

\medskip
Let $R>0,$ $y \in \Omega_{1,R}\cup \Gamma(R)$, $0<a<b$. Denote $A:=A_{a}^{b}$ and $A_\rho:= A_{a+\rho}^{b-\rho}$, similar to Lemma \ref{l.2.3a}, we have
\begin{lem}\label{l.2.3b} There exists a constant
$C>0$ depending only on $N$ such that for any  $u\in W^{3,m}_{0}(\Omega)$ and $0 < \rho < \min(1, \frac{b-a}{4}),$ we have
$$\|\D u\|_{L^m(A_\rho\cap \O)}^m\leq C \left(\frac{1}{\rho^{4}}\|u\|_{L^m(A\cap\O)}^m+\|\nabla(\D u)\|_{L^m(A\cap \O)}^m\right).$$
\end{lem}

\subsection{Explicit estimate via Morse index}

\medskip
  \begin{lem}\label{l.2.4a}
\label{lemnew3} Let $f$ satisfy $(H_1)$  and let $u \in W^{3,m}(\O)$ be  a weak solution to $(E_{3,m}),$  with Morse index $i(u) < \infty$. Then for any $y\in \Gamma(R)\cup\O_{1, R}$ with $R > 0$, there exists $j_{0}\in \{1,2,...,1+i(u)\}$ such that
\begin{align}
 \nonumber \int_{A_{j_0}\cap\Omega}|\nabla(\D u)|^m dx + \int_{A_{j_0}\cap\Omega}f(x,u)u dx \leq C
\left(\frac{1+i(u)}{R}\right)^\frac{3m(\mu+m)}{\mu}.
\end{align}
\end{lem}
 \noindent{\bf Proof.} Take $\eta \in C^6(\overline\O)$. By direct
 calculations, we get,
\begin{align*}
&\int_{\O} |\nabla(\D u)|^{m-2} |\nabla(\Delta (u\eta^{2k}))|^2 dx-(1+C\epsilon)\int_{\O} |\nabla(\Delta u)|^m\eta^{4k} dx\\
 &\leq \; C_\e
  \int_{\O} \Big[\rho^{-2m} |\nabla u|^m\eta^{4k-2m}+\rho^{-m}\left(|\D u|^m+|\nabla^{2} u|^{m}\right)\eta^{4k-m}
   + \rho^{-3m}|u|^{m} \eta^{4k-3m}\Big]dx.
\end{align*}

Using Lemma \ref{l.2.7a}, with  $\eta = \zeta$, we derive that
\begin{align*}
 \int_{\O} |\nabla(\D u)|^{m-2} |\nabla(\Delta (u\zeta^{2k}))|^2 dx
  \leq (1 + C\e)\int_{\O} |\nabla(\Delta u)|^m\zeta^{4k} dx + C_\e\rho^{-3m} \int_{\O} |u|^{m}\zeta^{4k-3m} dx.
\end{align*}
As in section 2, we can easily check that $\{u\phi_{j}^{m}\}_{1\leq j \leq i(u) + 1}$
are linearly independent, so there exists $j_{0}\in
\{1,2,...,1+i(u)\}$ such that $\Lambda_u(u \phi_{j_0}^{2k} )\geq 0$. The above estimate with $\zeta = \phi_{j_0}$ implies then
    \begin{align}\label{2.Da}
\displaystyle{\int_{\Omega}}f'
(x,u)u^{2}\phi_{j_{0}}^{4k}dx -\left(m-1\right)(1+C\epsilon)\int_{\Omega}|\nabla(\Delta u)|^m
\phi^{4k}_{j_{0}} dx \leq
\frac{C_\epsilon}{R^{3m}}(1+i(u))^{3m}\int_{\Omega}|u|^{m}\phi_{j_{0}}^{4k-3m} dx.
\end{align}
Now, take $u \phi_{j_0}^{4k}$ as the test function for $(E_{3,m})$, the integration by parts yields that
\begin{align*}
&\int_{\O}|\nabla(\Delta u)|^m \phi_{j_0}^{4k} dx - \int_{\O}f(x,u)u \phi_{j_0}^{4k} dx \\
 &\leq \; \int_{\O}|\nabla(\D u)|^{m-1}\left| \D u \nabla (\phi_{j_0}^{4k}) +2\nabla^{2} u \nabla (\phi_{j_0}^{4k})+\nabla u \D (\phi_{j_0}^{4k}) +
  2\nabla u \nabla^{2} (\phi_{j_0}^{4k})+ u \nabla(\D (\phi_{j_0}^{4k}))\right|.
\end{align*}
Applying Young’s inequality, we get, for any $\epsilon >0$
\begin{align*}
&(1 - C\epsilon)\int_{\O}|\nabla(\Delta u)|^m \phi_{j_0}^{4k} dx - \int_{\O}f(x,u)u \phi_{j_0}^{4k} dx\\
 &\leq \;C_\e
  \int_{\O} \Big[\rho^{-2m} |\nabla u|^m\phi_{j_0}^{4k-2m}+\rho^{-m}\left(|\D u|^m+|\nabla^{2} u|^{m}\right)\phi_{j_0}^{4k-m}
   + \rho^{-3m}|u|^{m} \phi_{j_0}^{4k-3m}\Big]dx .
\end{align*}

Using   Lemma \ref{l.2.7a}, with  $\eta = \phi_{j_0}$ we can conclude: For any $\epsilon>0$, there exists $C_\epsilon$
 such that
\begin{align}\label{2.La}
(1 - C\epsilon)\int_{\O}|\nabla(\Delta u)|^m \phi_{j_0}^{4k}dx-\int_{\O}f(x,u)u\phi_{j_0}^{4k}dx\leq
\frac{C_\epsilon}{R^{3m}}(1+i(u))^{3m}\int_{\Omega}|u|^{m}\phi_{j_{0}}^{4k-3m} dx.
\end{align}
Multiplying \eqref{2.La} by $\frac{(m-1)(1+ 2C\epsilon)}{1- C\epsilon}$
adding it with \eqref{2.Da}, we obtain  from $(H_1)$ that
  \begin{align*}
 &(m-1) C\epsilon\int_{\O}|\nabla(\Delta u)|^m \phi_{j_0}^{4k}dx+ \left(\mu-\frac{C\epsilon(1-2m)-m+2}{1-C\epsilon}\right)\int_{\O}f(x,u)u\phi_{j_0}^{4k}dx\\
 &\leq \;
\frac{C_\epsilon}{R^{3m}}(1+i(u))^{3m}\int_{\Omega}|u|^{m}\phi_{j_{0}}^{4k-3m} dx +C.
\end{align*}
Fix $\epsilon <\frac{\mu-m+2}{C(\mu+2m-1})$, we get
\begin{align*}
 \int_{\O}|\nabla(\Delta u)|^m \phi_{j_0}^{4k}dx+\int_{\O}f(x,u)u\phi_{j_0}^{4k}dx\leq
\frac{C_\epsilon}{R^{3m}}(1+i(u))^{3m}\int_{\Omega}|u|^{m}\phi_{j_{0}}^{4k-3m} dx + C.
\end{align*}

\medskip
 Now, we will proceed as the proof of \eqref{ecf2}, the claim follows. \qed

 \subsection{Proof of Theorem \ref{main2} for $k = 3$}
We show firstly the Pohozaev identity  for $(E_{3,m})$.
\begin{lem}\label{l.2.1a}
Let $u\in W^{3,m}(\O)$ be weak solution to $(E_{3,m})$. Let $\psi \in C_c^{4}(B_R(y))$. Then
\begin{align*}
& N\int_\O F(x,u)\psi dx +
 \int_\O \nabla_{x} F(x,u)\cdot n\psi dx
- \frac{N-3m}{m}\int_\O |\nabla(\Delta u)|^{m}\psi dx \\
 =&\; \frac{1-m}{m}\int_\O |\nabla(\Delta u)|^{m}(\nabla\psi\cdot n) dx  - \int_\O F(x,u)\nabla \psi\cdot n dx-2 \int_{ \Omega}\Delta u|\nabla(\Delta u)|^{m-2}\nabla(\Delta u)\cdot\nabla\psi dx\\
 - & \; \int_\O \D\psi\left[|\nabla(\Delta u)|^{m-2}\nabla(\Delta u)\right]\nabla(n\cdot\nabla u)dx-\int_{\Omega} \left[|\nabla(\Delta u)|^{m-2}\nabla(\Delta u)\cdot\nabla(\Delta\psi)\right](\nabla u\cdot n)dx\\
 -&\; 2 \int_\O \left[|\nabla(\Delta u)|^{m-2}\nabla(\Delta u)\right]\nabla\Big[\nabla^2u(n, \nabla\psi) + \nabla u\nabla\psi\Big]dx\\
+& \int_{\partial\Omega_R(y)}\left[|\nabla(\Delta u)|^{m-2}(\nabla(\Delta u)\cdot v)\cdot (\nabla (\Delta u)\cdot n)\psi-\frac{1}{m}(\nabla(\Delta u))^{m}( v\cdot n)\psi\right]d\sigma.
\end{align*}
\end{lem}

For the boundary terms, we have
\begin{lem}\label{l.2.2a}
 There exists $R_1 > 0$ depending only on $\O$ such that for any $u$ smooth function in $\in W^{3,m}_{0}(\O)$, any $0<R<R_1$, $y\in
\Gamma(R)$ and any nonnegative function $\psi$, there holds
\begin{align}
 \nonumber \int_{\partial\Omega_R(y)}\left[|\nabla(\Delta u)|^{m-2}(\nabla(\Delta u)\cdot v)\cdot (\nabla (\Delta u)\cdot n)\psi-\frac{1}{m}(\nabla(\Delta u))^{m}( v\cdot n)\psi\right] d\sigma \leq 0.
\end{align}
 \end{lem}

\noindent
\textbf{Proof.} Let  $m>2,$ and we proceed similarly as in the proofs of Lemma 2.2. in \cite{HHF,lec}, to show that there exists $R_1=R_1(\O)>0$ such that,
 if $0<R\leq R_1$ and $y \in \Gamma(R)$ then $v\cdot n \leq 0,$ for $x \in \partial\Omega_R(y).$
 If $\nabla u \neq 0$ for $x \in \partial\Omega_R(y)$, we have $v= \epsilon \frac{\nabla(\Delta u)}{|\nabla(\Delta u)|}$ with $\epsilon= \pm 1$.
  Therefore, there holds
 \begin{align*}
&\; \int_{\partial\Omega_R(y)}\left[|\nabla(\Delta u)|^{m-2}(\nabla(\Delta u)\cdot v)\cdot (\nabla (\Delta u)\cdot n)\psi-\frac{1}{m}(\nabla(\Delta u))^{m}( v\cdot n)\psi\right] d\sigma\\
 &=\; \left(1-\frac{1}{m}\right)\int_{\partial\Omega_R(y)}(\nabla(\Delta u))^{m}( v\cdot n)\psi\leq 0,\quad \forall\; x \in \p\O_R(y).
 \end{align*}
So we are done. \qed

\medskip
Similar to Proposition \ref{c.2.1}, we can claim
 \begin{prop}\label{c.2.1a}
There exists $R_0 > 0,$ $C > 0$ and $\rho> 0$   who satisfies the following property: Let $u\in W^{3,m}(\O)$
be a weak solution of $(E_{3,m})$ with $f$ satisfying $(H_{1})$--$(H_{3})$, let $0< R\leq R_0$, $y \in\Gamma(R)\cup \O_{1, R}$ and  $\psi \in C_c^6(B_R(y))$, $\psi \in [0, 1]$, there holds
\begin{align}
\label{xx.0.ve}
 \begin{split}
 &\; \int_\O f(x,u)u\psi dx + \int_\O |\nabla(\Delta u)|^2\psi dx\\
\leq&\; CR\rho^{-1}\int_{A_{R,\psi}(y)} f(x, u) udx
 + C\left(1+R\rho^{-1}+R^m\rho^{-m}\right)\int_{A_{R,\psi}(y)}|\nabla(\D u)|^m dx \\
 &+C\left(\rho^{-2m}+ CR^m\rho^{-3m}\right)\int_{A_{R,\psi}(y)}|\nabla u|^m dx+  C\left(\rho^{-m}+R^m\rho^{-2m}\right)\int_{A_{R,\psi}(y)}|\D u|^m dx \\
&+\;C\left(\rho^{-3m}+R^m\rho^{-4m}\right)\int_{A_{R,\psi}(y)}|u|^m dx+C\left(1+\rho^{-1}\right)R^{N+1}.
\end{split}
\end{align}
Here $C$ is a positive constant depending on $\O, N, k, \mu, \theta$; and $A_{R, \psi}(y) = B_{R}(y)\cap\Omega\cap\{\nabla \psi\neq0\}.$
\end{prop}
\textbf{Proof.} Using Lemmas \ref{l.2.1a}- \ref{l.2.2a},  $(H_{1})$--$(H_{3})$  and  by $\eqref{B}$, we obtain
\begin{align*}
\begin{split}
 &\; \frac{N-3m}{m}\left[(1+\theta)\int_{\O}f(x,u)u\psi dx - \int_{\O}|\nabla(\Delta u)|^m\psi dx \right]\\
\leq&\; CR\|\nabla\psi\|_\infty \int_{A_{R,\psi}(y)}\left[|\nabla(\Delta u)|^m + f(x, u) u \right]dx + CR\int_\O f(x, u)u\psi dx+C(\|\psi\|_\infty+ \|\nabla \psi\|_\infty)R^{N+1}\\
&+ \int_{A_{R,\psi}(y)} \Big|\left[|\nabla(\Delta u)|^{m-2}\nabla(\Delta u)\right]\nabla\Big[\nabla^2u(n, \nabla\psi)+CR\|\nabla(\Delta\psi)\|_\infty \int_{A_{R,\psi}(y)} |\nabla(\Delta u)|^{m-1}|\nabla u|dx \\
&\;+C\|\nabla\psi\|_\infty \int_{A_{R,\psi}(y)}|\Delta u||\nabla(\Delta u)|^{m-1} dx+ \nabla u\nabla\psi\Big]\Big| dx+ \int_\O \Big| \D\psi\left[|\nabla(\Delta u)|^{m-2}\nabla(\Delta u)\right]\nabla(n\cdot\nabla u)\Big| dx.
\end{split}
\end{align*}
 By Young's inequality, we get

\begin{align}
\label{xx.0.0Ajjva}
\begin{split}
&\frac{N-3m}{m}\left[(1+\theta)\int_{\O}f(x,u)u\psi dx - \int_{\O}|\nabla(\Delta u)|^m\psi dx \right]\\
&\leq \int_{A_{R,\psi}(y)} \Big||\nabla(\Delta u)|^{m-1}\nabla\Big[\nabla^2u(n, \nabla\psi)+\nabla u\nabla\psi\Big]\Big| dx+ \int_{A_{R,\psi}(y)} \Big| \D\psi\left[|\nabla(\Delta u)|^{m-1}\right]\nabla(n\cdot\nabla u)\Big|\\
&+CR\rho^{-1}\int_{A_{R,\psi}(y)}\left[|\nabla(\Delta u)|^m + f(x, u) u \right]dx + C\int_{A_{R,\psi}(y)}|\nabla(\D u)|^m dx+CR\int_\O f(x, u)u\psi dx\\
&+ CR^{m}\rho^{-3m}\int_{A_{R,\psi}(y)}|\nabla u|^m dx+ C\rho^{-m}\int_{A_{R,\psi}(y)}|\D u|^m dx+C(1+\rho^{-1} )R^{N+1}.
\end{split}
\end{align}

 \noindent
 We will use also the following lemma, which proof is given later.
\begin{lem}\label{l.2.7aa}
For any  $R < 1$, and  $\psi \in C_c^6(B_R(y))$, with $0 \leq\zeta \leq 1,$ there exists $C >0$ such that
\begin{align}
\label{xxx.013.0Ba}
\begin{split}
& \int_{A_{R,\psi}(y)} \Big||\nabla(\Delta u)|^{m-1}\nabla\Big[\nabla^2u(n, \nabla\psi) + \nabla u\nabla\psi\Big]\Big| dx + \int_{A_{R,\psi}(y)} \Big| \D\psi\left[|\nabla(\Delta u)|^{m-1}\right]\nabla(n\cdot\nabla u)\Big| dx\\
 \leq &C\left(1+R^m\rho^{-m}\right)\int_{A_{R,\psi}(y)}|\nabla(\D u)|^m dx + C\left(\rho^{-m}+R^m\rho^{-2m}\right)\int_{A_{R,\psi}(y)}|\D u|^m dx \\
 &+\left(\rho^{-2m}+ CR^m\rho^{-3m}\right)\int_{A_{R,\psi}(y)}|\nabla u|^m dx + C\left(\rho^{-3m}+R^m\rho^{-4m}\right)\int_{A_{R,\psi}(y)}|u|^m dx  .
\end{split}
\end{align}
\end{lem}

Combining \eqref{xx.0.0Ajjva}-\eqref{xxx.013.0Ba}, there holds

\begin{align}
\label{xx.0.0Afjjva}
\begin{split}
 &\; (1+\theta)\int_{\O}f(x,u)u\psi dx - \int_{\O}|\nabla(\Delta u)|^m\psi dx\\
\leq&\; CR\|\nabla\psi\|_\infty \int_{A_{R,\psi}(y)}\left[|\nabla(\Delta u)|^m + f(x, u) u \right]dx + CR\int_\O f(x, u)u\psi dx\\
 \leq & \;C\left(1+R^m\rho^{-m}\right)\int_{A_{R,\psi}(y)}|\nabla(\D u)|^m dx +  C\left(\rho^{-m}+R^m\rho^{-2m}\right)\int_{A_{R,\psi}(y)}|\D u|^m dx \\
 & \; +\left(\rho^{-2m}+ CR^m\rho^{-3m}\right)\int_{A_{R,\psi}(y)}|\nabla u|^m dx + C\left(\rho^{-3m}+R^m\rho^{-4m}\right)\int_{A_{R,\psi}(y)}|u|^m dx  .
\end{split}
\end{align}

\medskip
 Take $u \zeta^{4k}$ as the test function for $(E_{3,m})$, using Lemmas \ref{l.2.7a}, for any $\epsilon>0$ there exists $C_\epsilon$ such that
\begin{align}\label{xx.0.0Aa}
\int_{\O}|\nabla(\Delta u)|^m \zeta^{4k}dx-\int_{\O}f(x,u)u\zeta^{4k}dx\leq
C\e\int_{\O}|\nabla(\Delta u)|^m\zeta^{4k} dx + C_\e R^{-3m}\int_{\O}|u|^m dx.
\end{align}
Let $0<\theta'<\theta,$ and $\psi = \zeta^{4k}$. choosing $\e, R > 0$ small enough, and  Combining \eqref{xx.0.0Afjjva}-\eqref{xx.0.0Aa}, we have \eqref{xx.0.ve}. \qed

  \bigskip
\noindent{\bf Proof of Lemma \ref{l.2.7aa}.}   Firstly we have ,

\begin{align*}
\int_{A_{R,\psi}(y)}|\nabla^3(u\nabla\psi)|^m dx \leq &\; C \int_{{A_{R,\psi}(y)}}|\nabla\D(u\nabla\psi)|^mdx \\
&\leq C\int_{A_{R,\psi}(y)}\Big(\left(|\nabla^2 u|^m +|\D u|^m\right)|\nabla^2\psi|^m+ |\nabla u|^m|\nabla^3\psi|^m + |u|^m|\nabla^4\psi|^m\Big)dx\\
& + C\int_{A_{R,\psi}(y)}|\nabla(\D u)|^m|\nabla\psi|^m dx.
\end{align*}

Using the inequality $$|\nabla^3u \nabla\psi|\leq |\nabla^3(u\nabla\psi)| +3|\nabla^2u||\nabla^2\psi|+3|\nabla u||\nabla^3\psi|+|u||\nabla^4\psi|,$$  we get
 \begin{align*}
&\Big|\left[|\nabla(\Delta u)|^{m-1}\right]\nabla\Big[\nabla^2u(n, \nabla\psi) + \nabla u\nabla\psi\Big]\Big|+\Big| \D\psi\left[|\nabla(\Delta u)|^{m-1}\right]\nabla(n\cdot\nabla u)\Big|\\
 \leq & \; CR|\nabla(\D u)|^{m-1}\Big(|\nabla^3(u\nabla\psi)| + |\nabla^2u||\nabla^2\psi| + |\nabla u||\nabla^3\psi| + |u||\nabla^4\psi|\Big)\\
  &\;+ C|\nabla(\D u)|^{m-1}\Big(|\nabla^2 u||\nabla\psi| + |\nabla u||\nabla^2\psi|\Big) +C|\nabla(\D u)|^{m-1}\Big[R|\nabla^2 u||\Delta \psi| + |\nabla u||\Delta\psi||\Big].
\end{align*}
Combining all these inequalities,  using $|\D\psi|\leq|\nabla^2\psi|,$ and by Young's inequality,  we arrive at
\begin{align*}
\begin{split}
&  \int_{A_{R,\psi}(y)} \Big||\nabla(\Delta u)|^{m-1}\nabla\Big[\nabla^2u(n, \nabla\psi) + \nabla u\nabla\psi\Big]\Big| dx + \int_{A_{R,\psi}(y)} \Big| \D\psi\left[|\nabla(\Delta u)|^{m-1}\right]\nabla(n\cdot\nabla u)\Big| dx \\
 \leq & \;C\left(1+\rho^{-m}\right)\int_{A_{R,\psi}(y)}|\nabla(\D u)|^m dx +  C\int_{A_{R,\psi}(y)}|\nabla^2u|^m\left(R^m|\nabla^2\psi|^m+|\nabla\psi|^m\right) dx \\
 & \; + C\left(\rho^{-2m}+R^m\rho^{-3m}\right)\int_{A_{R,\psi}(y)}|\nabla u|^m dx + CR^m\rho^{-4m}\int_{A_{R,\psi}(y)}|u|^m dx + CR^m\rho^{-2m}\int_{A_{R,\psi}(y)}|\D u|^m dx.
\end{split}
\end{align*}
By direct calculations, we get,
\begin{align*}
|\nabla^2u|^{m}|\nabla^{n} \psi|^{m}
 \leq & \; C_{m}\Big[ |\nabla^2(u\nabla^{n} \psi)|^{m}+|\nabla u|^m|\nabla^{n+1}\psi|^m+|u|^m|\nabla^{n+2}\psi|^m\Big].
\end{align*}
As $u\nabla^{2}\psi=0$ on $\partial\O,$ there exists $C>0$ depending only on $\O$ such that
\begin{align*}
\int_{A_{R,\psi}(y)}|\nabla^2(u\nabla^{n} \psi)|^{m}  \leq  C\int_{{A_{R,\psi}(y)}} |\D(u\nabla^{n} \psi)|^{m}dx.
\end{align*}
Combining the last tow inequality, there holds
\begin{align*}
\begin{split}
 & \int_{A_{R,\psi}(y)}|\nabla^2u|^m\left(R^m|\nabla^2\psi|^m+|\nabla\psi|^m\right) dx \\
\leq & \;C\left(\rho^{-m}+R^m\rho^{-2m}\right)\int_{A_{R,\psi}(y)}|\D u|^m dx +C\left(\rho^{-2m}+R^m\rho^{-3m}\right)\int_{A_{R,\psi}(y)}|\nabla u|^m dx\\
& \;+C\left(\rho^{-3m}+R^m\rho^{-4m}\right)\int_{A_{R,\psi}(y)}|u|^m dx  .
\end{split}
\end{align*}
Combining all these inequalities,  the estimate \eqref{xxx.013.0Ba} is proved. \qed

\bigskip
\noindent
{\bf Proof of Theorem \ref{main2} for $k=3$ completed}.
 Now, we are in position to prove Theorem \ref{main2} for $k=3$ . Fix $$R = R_0, \quad 4k = 3m+\frac{3m^{2}}{\mu}, \quad \rho:=\frac{R}{10(i(u)+1)}, \quad A_{j_0,\rho}:=A_{a_{j_0}+\rho}^{b_{j_0}-\rho} \subset A_{j_0} \;\mbox{be as in $(*)$}.$$
Using Remark \ref{rem2.1} and  lemma \ref{l.2.4a}, there holds
\begin{eqnarray}\label{PPLj}
 \|u\|_{L^m(A_{j_{0}}\cap\O)}^m \leq C\left(\int_{A_{j_{0}}\cap \O}f(x,u)u\right)^{\frac{m}{m+\mu}} + C
 \leq C(1+i(u))^{\frac{3m^{2}}{\mu}}.
 \end{eqnarray}
According to Lemmas \ref{l.2.3a}, \ref{l.2.3b}, \ref{l.2.4a} and \eqref{PPLj},  we can claim
\begin{align}\label{3.7b}
 \|\nabla(\D u)\|^m_{L^m(A_{j_0,\rho}\cap\O)}+ \|\nabla u\|^m_{L^2(A_{j_0,\rho}\cap\O)}+\|\ \D u\|^m_{L^m(A_{j_0,\rho}\cap\O)}\leq C
(1+i(u))^\frac{3m(\mu+m)}{\mu}.
\end{align}
Combining \eqref{xx.0.ve}, \eqref{PPLj} and  \eqref{3.7b}, one obtains
\begin{align*}
 \int_{\O}f(x,u)u\xi_{j_0} dx + \int_{\O}|\nabla(\Delta u)|^m \xi_{j_0} dx
\leq C(1+i(u))^\frac{3m(2\mu+m)}{\mu}.
\end{align*}
As $\frac{R}{2}< a_{j_{0}}$ , we get then for any $y \in \Gamma(R)\cup \O_{1, R}$,
\begin{align*}
\int_{B_{\frac{R_0}{2}}(y)\cap \O} \Big[|\nabla(\Delta u)|^m  + f(x,u)u\Big] dx \leq C
(1+i(u))^\frac{3m(2\mu+m)}{\mu}.
\end{align*}
The proof is completed by the covering argument.\qed

\end{document}